%% file: main-FINAL.tex
\documentclass[12pt,reqno]{amsart}

\usepackage{tikz}
\usepackage[T1]{fontenc}
\usepackage{lmodern}
\usepackage{microtype}
\usepackage{needspace}
\usepackage{a4wide}
\usepackage{amsmath,amsthm,amsfonts,amssymb,amscd}
\usepackage{enumitem}
\usepackage{graphicx}
\usepackage{xcolor}
\usepackage{float}
\numberwithin{equation}{section}
\usepackage[
  pagebackref,
  colorlinks=true,
  linkcolor=blue,
  citecolor=black,
  urlcolor=blue
]{hyperref}
\hypersetup{
  pdftitle={Tetrahedral points and global rigidity for Almgren-minimizing sets},
  pdfauthor={} ,
  colorlinks=true,
    citecolor=red
}

\usepackage{caption}
\theoremstyle{plain}
\newtheorem{theorem}{Theorem}[section]
\newtheorem{lemma}[theorem]{Lemma}
\newtheorem{proposition}[theorem]{Proposition}
\newtheorem{corollary}[theorem]{Corollary}
\theoremstyle{definition}
\newtheorem{definition}[theorem]{Definition}
\newtheorem{remark}[theorem]{Remark}

\newcommand{\Sp}{\mathbb{S}}

\newcommand{\R}{\mathbb{R}}

\newcommand{\Z}{\mathbb{Z}}

\newcommand{\ha}{\mathcal{H}}

\newcommand{\id}{{\rm id}}
\newcommand{\cH}{\mathcal{H}}

\newcommand{\interior}{\operatorname{int}}
\newcommand{\Lip}{\operatorname{Lip}}

\DeclareMathOperator{\dist}{dist}

\makeatletter
\def\l@subsection{\@tocline{2}{0pt}{2.5pc}{5pc}{}}
\def\l@subsubsection{\@tocline{2}{0pt}{5pc}{7.5pc}{}}
\makeatother

\author[M. Caselli]{Michele Caselli}
\author[Z. Lihn]{Zachary Lihn}

\address[M. Caselli]{Princeton University --  Fine Hall, 304 Washington Rd, Princeton, NJ 08540 
	\newline\indent 
	University of Sydney -- Quadrangle A14, Camperdown, NSW 2006, Australia}
\email{\href{mailto:mc3147@princeton.edu}{mc3147@princeton.edu}}

\address[Z. Lihn]{Princeton University -- Fine Hall, 304 Washington Rd, Princeton, NJ 08540 }
\email{\href{mailto:zl3140@princeton.edu}{zl3140@princeton.edu}}

\title{Tetrahedral points and classification of Almgren-minimizing sets in $\mathbb{R}^3$}

\input{y_mesh_support.tex}

\begin{document}
\begin{abstract}
We prove that a two-dimensional Almgren-minimizing set in three dimensions which is sufficiently close in Hausdorff distance to a tetrahedral cone contains an interior tetrahedral point. As a consequence, we obtain quantitative regularity near tetrahedral cones and we complete the classification of entire Almgren-minimizers in $\mathbb{R}^3$. 
\end{abstract}
\maketitle
\section{Introduction}\label{sec:introduction}

Two-dimensional Almgren-minimizing sets provide a variational model for soap films. A celebrated theorem by Taylor \cite{Almgren1976,Taylor1976}, see also \cite{David2009,David2010}, shows that these sets are locally $C^{1,\alpha}$-diffeomorphic to a plane $P$, a triple junction $Y$, or a tetrahedral cone $T$, and that these three models are the only Almgren-minimizing cones in three dimensions. Taylor's theorem is a local statement in that the description holds at a sufficiently small scale around each point, and this scale may depend on the point. 

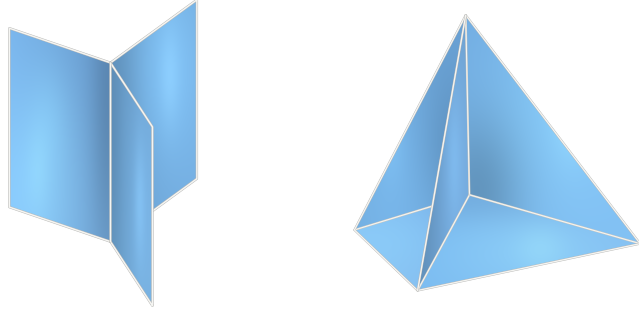
\begin{figure}[ht]
    \centering
     \resizebox{0.3\textwidth}{!}{\input{y_set.tikz}}
     \hspace{0.02\textwidth}
      \resizebox{0.27\textwidth}{!}{\input{t_set.tikz}}
    \caption{A triple junction $Y$ and a  tetrahedral cone \(T\).}
    \label{fig: Y and T}
\end{figure}

A full $\varepsilon$-regularity theorem asks for a stronger property: being close to one of these cones $\{P, Y, T\}$ in Hausdorff distance at one scale should imply that the set is a $C^{1,\alpha}$ version of the same cone at a fixed smaller scale. Such a property is known near planes and triple junctions \cite{David2009,David2010}. In these cases, closeness to the cone in Hausdorff distance guarantees the existence of a point of the same type in a smaller ball, and the regularity theory then gives the desired quantitative description.

For tetrahedral cones, the missing ingredient is precisely this existence property. Regularity near a point already known to be $T$ follows from Taylor's theorem; what remained open was to obtain the existence of such a point from Hausdorff closeness alone, without any a priori topological assumption on the set; see \cite[Remark 4.78]{David2013TaylorBeginners}.

In this work we prove that an Almgren-minimizing set sufficiently close to a tetrahedral cone in a ball contains a tetrahedral point in a smaller ball. This problem is explicitly stated in several places in the literature: \cite[Remark 4.78]{David2013TaylorBeginners}, \cite[Section 4.c]{DavidLecNotes}, \cite[Remark 18.43]{David2009}, \cite[Remarks 8.4 and 8.5]{ColomboEdelenSpolaor2022}, and the discussion below \cite[eq. (1.11)]{Liang2016GlobalPlanes}. Our main theorem provides the missing ingredient and yields an $\varepsilon$-regularity theorem near tetrahedral cones.

The existence of $T$ points of Theorem \ref{thm:existence} also completes the classification of entire two-dimensional Almgren-minimizing sets in $\mathbb R^3$. Blow-down limits of such sets are minimal cones, and the cases of planar and triple-junction blow-downs have already been classified; see \cite{Liang2014,David2013TaylorBeginners}. In the remaining case of a tetrahedral blow-down, our existence theorem gives a point whose density agrees with that of the blow-down. The monotonicity formula then forces the set itself to be a tetrahedral cone.
 
 We now state our main existence theorem. Let
\begin{equation*}
    d_{x,r}(E,F)=\frac{1}{r}\max\left\{
    \sup_{z\in E\cap B_{r}(x)}\operatorname{dist}(z,F),\,
    \sup_{z\in F\cap B_{r}(x)}\operatorname{dist}(z,E)\right\} 
\end{equation*}
be the normalized (local) Hausdorff distance between closed sets.

\begin{theorem}\label{thm:existence}
There exists a universal $\varepsilon_0>0$ with the following property. Let $U\subset\mathbb R^3$ be open, $ x \in U$ and $r>0$ such that $\overline B_{r}(x)\subset U$. Let $E\subset U$ be a closed, reduced, two-dimensional Almgren minimal set. Assume there exists a tetrahedral cone $T_x$ centered at $x$ such that
\begin{equation*}
d_{x,r}(E,T_x)<\varepsilon_0.
\end{equation*}
Then $E\cap B_{r/4}(x)$ contains a $T$ point. 
\end{theorem}

    Note that we do not assume any topological condition on $E$ (cf. \cite[Proposition 18.1]{David2009}) nor the existence of a singular point; in particular, $x$ need not belong to $E$.

\subsubsection*{Classification of entire minimizers.}

 Theorem~\ref{thm:existence} and the work of \cite{Taylor1976,David2009,David2010} immediately imply the classification of global Almgren-minimizing sets in $\mathbb R^3$. This fact is known for Mumford-Shah minimal sets \cite[Theorem 1.9]{David2009} and for size-minimizing currents \cite[Theorem 4.1]{MorganInvent}, and the only case that remained open for Almgren-minimal sets is that of sets with blow-down equal to a tetrahedral cone. 

 Indeed, consider
 \[
\theta_E(p,0):=\lim_{r \downarrow0}\theta_E(p,r), \qquad  
\theta_E(\infty) :=\lim_{r\to\infty}\theta_E(p,r), \qquad \theta_T :=\cH^2(T\cap B_1) , 
\]
where $\theta_E(p,r):=r^{-2}\cH^2(E\cap B_r(p))$ is the density ratio in the monotonicity formula, and suppose that an entire minimal set $E \subset \R^3$ has a tetrahedral blow-down. Applying Theorem~\ref{thm:existence} in a large ball we get the existence of a tetrahedral point $p_T \in E$. Since the densities of $E$ at $p_T$ and at infinity are both equal to $\theta_T$, the monotonicity formula gives
\[
    \theta_T = \theta_E(p_T,0)
    \leq \theta_E(p_T,r)
    \leq \theta_E(\infty) = \theta_T , 
    \qquad \text{for every } r>0.
\]
Hence, the density ratio is constant, and the equality case of the monotonicity
formula implies that $E$ is a tetrahedral cone. Hence, our result completes the classification of entire Almgren-minimal sets in $\R^3$.

 \begin{theorem}\label{thm:classification}
Let $E\subset\mathbb R^3$ be a closed, reduced, two-dimensional Almgren minimal set in $\R^3$. Then, $E$ is either a plane, a triple-junction cone, or a tetrahedral cone.
\end{theorem}

\subsubsection*{Regularity near tetrahedral cones}

Theorem \(\ref{thm:existence}\) also supplies the missing hypothesis in the existing regularity theory near \(T\). Once a tetrahedral point has been found, the estimates in \cite{Taylor1976,David2009,David2010} yield a \(C^{1,\alpha}\) parametrization on a smaller ball. 
 
\begin{corollary}\label{cor:epsilon}
Let $U \subset \R^3$ be open with $0\in U$ and $\overline B_{2r} \subset U$, and let $E\subset U$ be a reduced two-dimensional Almgren-minimal set in $U$. Let $T$ be the  tetrahedral cone centered at the origin. There exist universal $\varepsilon_0>0$ and $\alpha\in(0,1)$ such that if
\begin{equation*}
    d_{0,2r}(E,T)<\varepsilon_0,
\end{equation*}
then there is a \(C^{1,\alpha}\) diffeomorphism
\( \Phi:B_{2r} \to\Phi(B_{2r}) \)
satisfying 
\(E\cap B_r=\Phi(T\cap B_{2r})\cap B_r\) and $ \|\Phi-\operatorname{id}\|_{L^\infty(B_{2r})}\lesssim r .$
\end{corollary}

\begin{remark}
    Theorem \ref{thm:existence} extends to almost-minimal sets with gauge function $h(t)\le Ct^\beta$, $\beta>0$, provided the gauge is sufficiently small at the scale under consideration. Precisely, for every $\beta>0$ there exists $\varepsilon_\beta>0$ such that, under the remaining hypotheses of Theorem \ref{thm:existence}, $d_{x,r}(E,T_x)+Cr^\beta<\varepsilon_\beta$ implies that $E\cap B_{r/4}(x)$ contains a $T$-point. In the normalized contradiction argument, the additional almost-minimality error in the final area comparison vanishes. For clarity of presentation, we state and prove the theorem only for minimal sets.
\end{remark}

\subsubsection*{The history of possible counterexamples.} 
\begin{figure}[ht]
    \centering
    \includegraphics[width=0.3\textwidth]{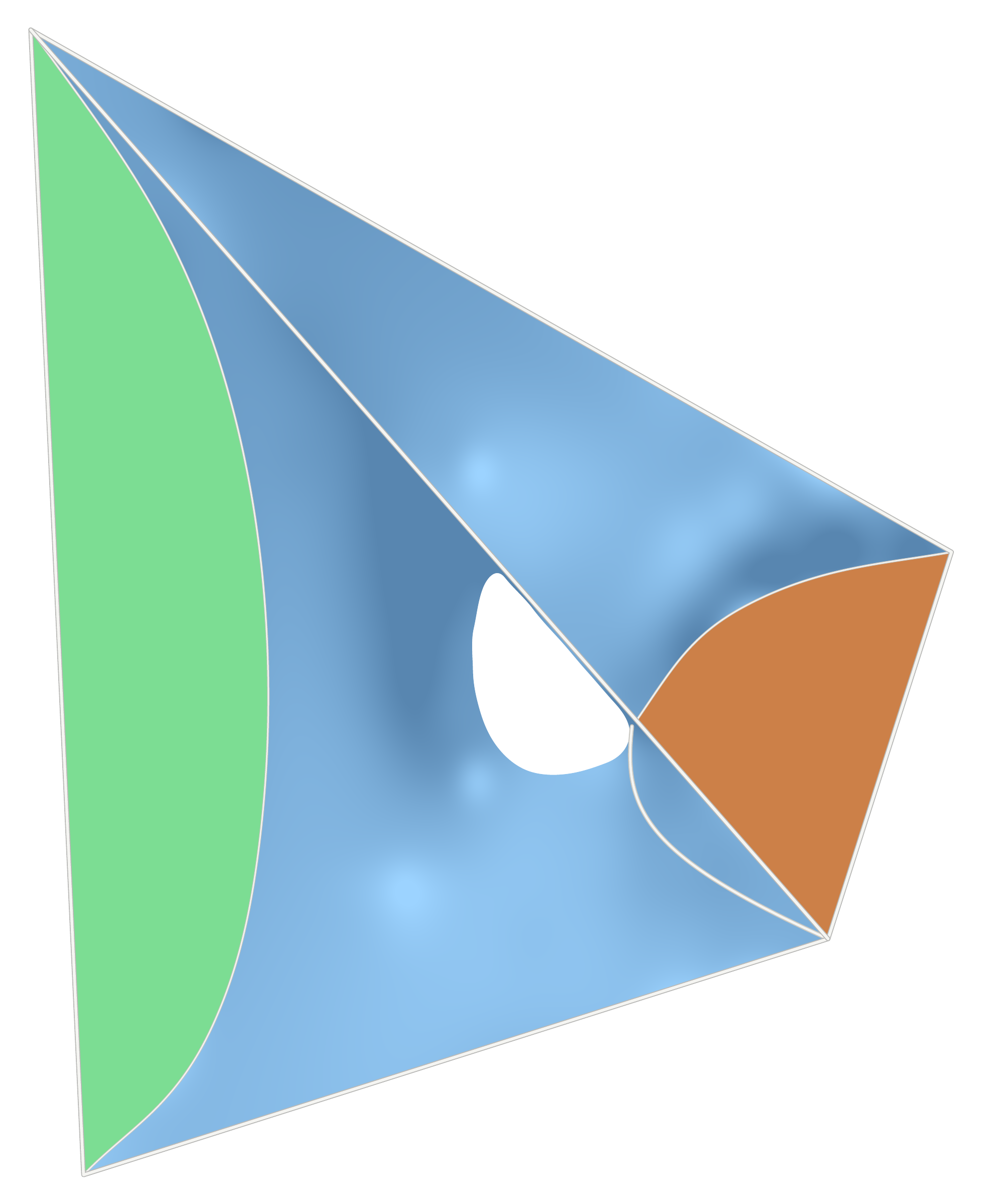}
\hspace{0.02\textwidth}
\includegraphics[width=0.3\textwidth]{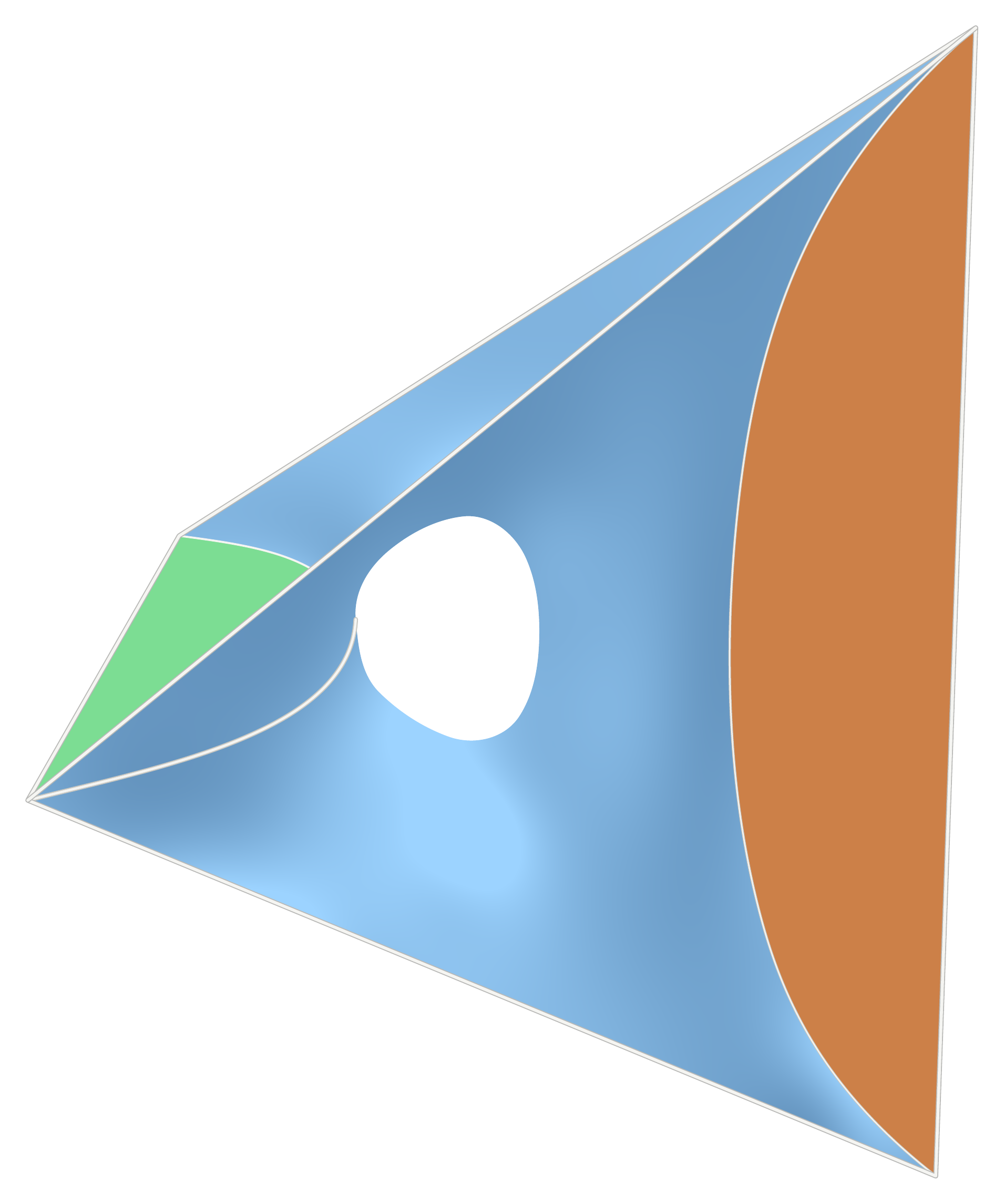}
    \caption{A set spanning the $1$-skeleton of the tetrahedron and with singularities only of type $Y$.}
    \label{fig:first top example}
\end{figure}

The existence of a tetrahedral point does not follow from the boundary topology and the known local structure alone. Examples of sets that coincide with $T$ on the $1$-skeleton of the tetrahedron and have singularities only of type $Y$ have been known for many years \cite{LawlorMorgan1994PairedCalibrations,Morgan1995GeometricMeasureTheory,HuffR,BellettiniPaoliniPasquarelli2018TripleCovers}; see Figure \ref{fig:first top example}.

\begin{figure}[H]
    \centering
    \includegraphics[width=0.8\textwidth]{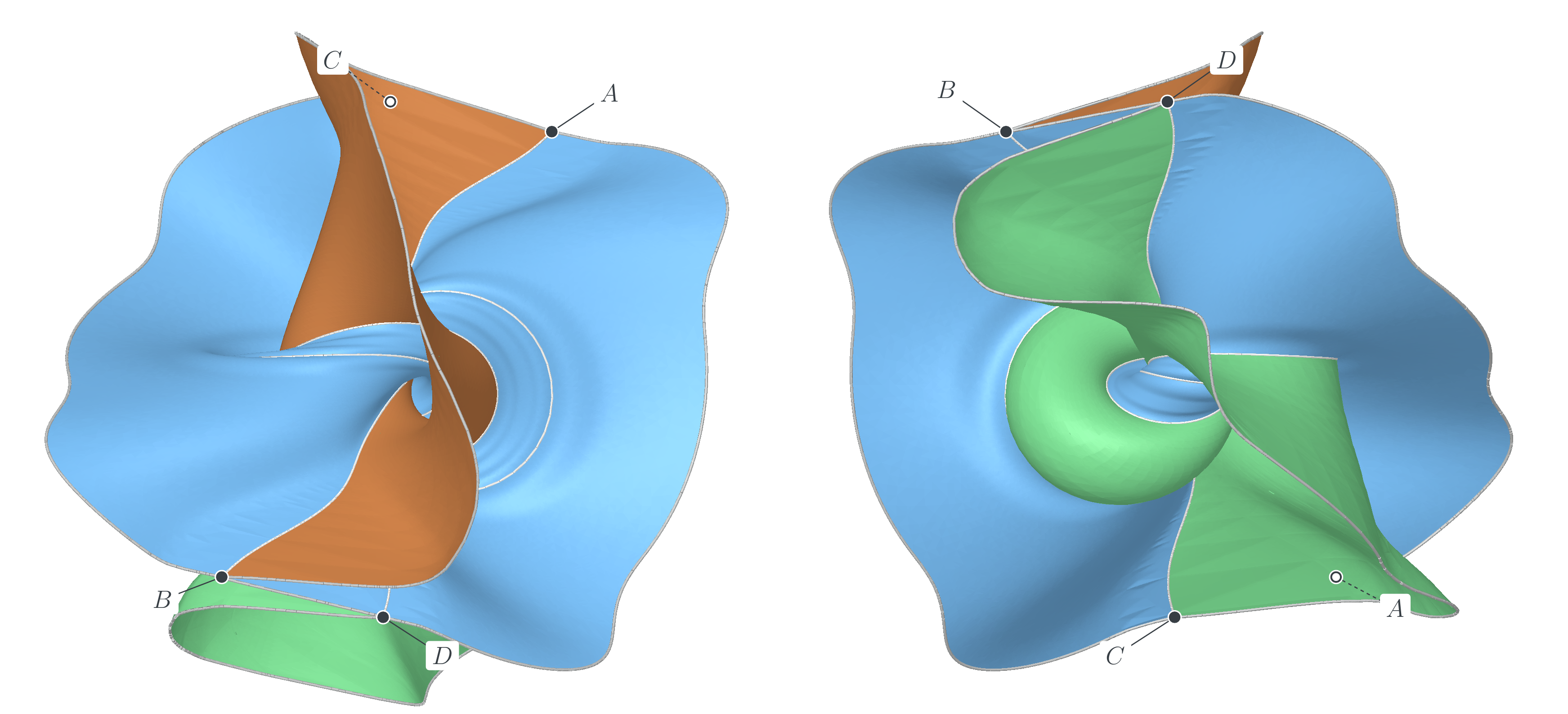}
    \caption{A set whose boundary is topologically the $1$-skeleton of a tetrahedron, with singularities only of type $Y$, and with knotted singular curves. Here $A,B,C,D$ are the vertices of the tetrahedron and three sheets meet at each of these vertices.}
    \label{fig: Liang knotted}
\end{figure}

This type of possible counterexample was later removed in \cite{Liang2014}, where the author shows that an Almgren-minimizing set with such a property must have the resulting singular curves knotted with respect to their boundary. These surfaces are harder to visualize (attempt in Figure \ref{fig: Liang knotted}) and, as already pointed out in \cite{Liang2014}, it was believed to be unlikely that an Almgren-minimizing set with such a topology exists. Indeed, by our Theorem \ref{thm:existence} a set with this topology cannot be Almgren-minimizing. 

\subsection{Outline of the proof}

We argue by contradiction, following an idea in \cite{Liang2014}.
The absence of tetrahedral points forces all interior singularities of $E$ to be of type $Y$. Hence there is a regular curve in the complement of $E$ joining two distinct faces of the tetrahedron, at a positive distance
from $E$. 

If the curve is unknotted, the argument in \cite{Liang2014} allows one to construct a deformation sending the portion of $E$ inside the tetrahedron into two boundary faces (together with the remaining edges of the tetrahedral $1$-skeleton). This target has a strictly smaller area than the $T$ cone inside the tetrahedron, yielding a contradiction since Hausdorff convergence implies the convergence of the masses. 

Nevertheless, we will show in Proposition~\ref{prop:deformation-key} that for general and simple topological reasons, the knottedness of the curve poses no obstacle to the existence of a Lipschitz map sending the portion of $E$ inside the tetrahedron into the same boundary target as in \cite{Liang2014}; see also \cite[Remark 8.5]{ColomboEdelenSpolaor2022}. Hence, we can conclude without any additional hypothesis.

\begin{remark}
Our technique also provides a conditional result on the uniqueness problem with prescribed tetrahedral boundary data; see e.g. \cite[Appendix D]{MaggiNovackRestrepo2026}. This is done in Theorem \ref{thm: prescribed-boundary} in the Appendix and assumes in addition that the minimal set is a $C^1$ version of $T$ in a neighborhood of the prescribed boundary $T\cap\partial Q$, namely the $1$-skeleton of the tetrahedron. This assumption fixes the topology near the boundary. To remove it, one should know a precise structure theorem for Almgren-minimizing sets near the boundary, which seems to be a very delicate problem \cite{CamilloFederico26}. 
\end{remark}

\subsubsection*{AI usage statement} The authors used OpenAI’s GPT-6 Astra in developing this work. The algebraic construction for the existence of $r:M \to A$ in the proof of Proposition \ref{prop:deformation-key} was suggested by GPT-6 Astra, while the geometric construction was found by the authors. The authors assume full responsibility for the correctness of the work. Access to the model was provided through OpenAI’s ChatGPT for Academic Researchers program.

\section{Preliminaries}\label{sec:preliminaries}

We begin by reviewing the definition of Almgren-minimal set and the results of \cite{Taylor1976,David2009,David2010} which will be necessary in the proof. For a friendly introduction to Almgren-minimal sets and Taylor's theorem we refer to \cite{David2013TaylorBeginners}.

\begin{definition}[Almgren-minimal set]\label{def:minimality}
Let $U\subset\mathbb R^3$ be open. We say that a relatively closed set $E\subset U$ is Almgren minimal if $E$ has locally finite $\cH^2$ measure in $U$ and the inequality
\begin{equation*}
    \cH^2(E\cap V)\leq \cH^2(\varphi(E)\cap V)
\end{equation*}
holds for every open set $V\Subset U$ and every Lipschitz map $\varphi : U \to U $ with $ \overline{W \cup \varphi(W)} \Subset V$, where $ W:= \{x  \mid  \varphi(x) \neq x \}$. 

We also say that $E$ is \emph{reduced} if 
\begin{equation*}
    \ha^2(E \cap B_{r}(x)) >0, \quad \mbox{for every } x\in E \mbox{ and } r>0. 
\end{equation*}
\end{definition}

\begin{remark} This is a minor variation of Almgren's original definition, also used in \cite{Liang2014, CamilloFederico26}, and it is slightly different from \cite[Definition 4.8 and Proposition 4.10]{David2009} which requires the existence of a path of maps $\{\varphi_t\}_{t\in [0,1]}$ starting from the identity to the final deformation. In a convex domain, any compactly supported final Lipschitz map induces such a family since its straight-line homotopy remains within that domain. In a general domain, this homotopy may leave the domain. Our proof also applies to the homotopy-based definition: the competitors constructed below admit homotopies supported in a ball compactly contained in the domain.

\end{remark}

\subsubsection*{The tetrahedral cone.} We now recall the geometry of the tetrahedral cone.
Let $v_1,v_2,v_3,v_4$ be the vertices of the regular tetrahedron centered at the origin with circumradius one. Thus,
\begin{equation*}
    |v_i| =1, \quad v_i \cdot v_j = -\frac{1}{3} \; (i\neq j), \quad v_1+v_2+v_3+v_4=0.
\end{equation*}
Let 
\begin{equation*}
    Q = \mathrm{conv} \{v_1,v_2,v_3,v_4\}
\end{equation*}
be the solid tetrahedron, and let $F_i = \mathrm{conv} \{v_j : j\neq i\}$ be the closed triangular face of $\partial Q$ \emph{opposite} to $v_i$. For $i \neq j$, we let
\begin{equation*}
    G_{ij} = \partial Q \setminus\bigl( \operatorname{relint}(F_{i}) \cup \operatorname{relint}(F_{j}) \bigr) . 
\end{equation*}
The Almgren-minimizing cone $T$ is the cone from the origin over the six edges of $Q$:
\begin{equation*}
    T = \bigcup_{1\leq i < j\leq 4} \left\{sv_i+ tv_j : s,t\geq 0 \right\}.
\end{equation*}
The cone $T$ has six faces. Its $Y$-point set, denoted by $T_Y$,
consists of the four singular rays with the origin removed; the origin
is the unique $T$ point. The rays point towards the vertices
$\{v_1,\ldots,v_4\}$, and
\begin{equation*}
     T_Y\cap Q=\bigcup_{i=1}^{4}\{t v_i:0<t\le1\}.
\end{equation*}

The six conical triangles inside $Q$ each have area $\sqrt{2}/3$. Therefore
\begin{equation*}
    \cH^2(T\cap Q)=2\sqrt{2}.
\end{equation*}
On the other hand each face $F_i$ has area $\cH^2(F_i) =2\sqrt{3}/3$. Therefore
\begin{equation*}
    \cH^2 (G_{ij}) = \frac{4\sqrt{3}}{3} <2\sqrt{2} = \cH^2(T\cap Q) . 
\end{equation*}
This is the strict inequality that we will exploit in the proof.


\subsubsection*{Regularity of Almgren-minimizing sets.}
The next lemma records the main consequences of \cite{David2009,David2010,DavidDePauwToro2008} we will need. It is a consequence of the $C^{1,\alpha}$ regularity of Almgren-minimizing sets near planes and triple junctions proved in \cite{Taylor1976,David2009,David2010}, along with a gluing procedure for such diffeomorphisms as described in \cite[Section 10]{DavidDePauwToro2008}, for example. The measure convergence \eqref{eq: meas conv} uses the lower and upper
semicontinuity results in \cite[Lemmas 3.3 and 3.12]{David2009}.
\begin{lemma} \label{lem:annular-diffeo}
    Let \(E_k\subset B_4 \) be closed, reduced, two-dimensional Almgren-minimal sets with $E_k \to T$ locally in Hausdorff distance in \(B_4\).
    There exists \(\alpha\in(0,1)\) such that, for all sufficiently large \(k\), there exists an ambient \(C^{1,\alpha}\) diffeomorphism \(\Phi_k:\mathbb R^3\to\mathbb R^3\) satisfying
    \begin{gather*}
         \Phi_k=\id\quad \text{on } B_{1/8} \cup\bigl(\mathbb R^3\setminus B_2 \bigr) ,  \\
         \|\Phi_k-\id\|_{L^\infty(\mathbb R^3)}+\|D\Phi_k-I\|_{L^\infty(\mathbb R^3)} \to 0 , \\
         \Phi_k^{-1}(E_k)\cap\left(  B_{3/2} \setminus\overline{B_{1/4}}\right)=T\cap\left( B_{3/2}\setminus\overline{B_{1/4}}\right) . 
    \end{gather*}
    Set
    $ X_k:=\Phi_k^{-1}(E_k)\cap Q$ and $Q_k:=\Phi_k(Q).$
    Then \(X_k=T\cap Q\) in a  neighborhood of \(\partial Q\) and
    \begin{equation}\label{eq: meas conv}
        \mathcal H^2(X_k)\to \mathcal H^2(T\cap Q)=2\sqrt2, \qquad \mathcal H^2(E_k\cap Q_k)\to 2 \sqrt2.
    \end{equation}
    Moreover, for some $C>0$ independent of $k$, 
    \begin{equation*}
        \dist(z,T\cap\partial Q)\le C \dist(z,\partial Q)  , \quad \mbox{for every } z\in X_k. 
    \end{equation*}
\end{lemma}

\section{Proof of Theorem~\ref{thm:existence}}

We begin by recalling the argument of \cite[Proposition 18.1]{David2009} and \cite[pp. 210-211]{Liang2014}, which says that the nonexistence of $T$-points forces the existence of a curve joining two faces of $\partial Q$.
\begin{lemma} \label{lem:topological-path}
     Let \(X\subset Q\) be closed set which agrees with \(T\cap Q\) in a  neighborhood of \(\partial Q\). Suppose that, at every point of \(X\cap\operatorname{int}(Q)\), the set \(X\) is locally diffeomorphic to a plane or a triple junction. Then there exist \( i \neq j \) and a \(C^1\) embedded curve
    \begin{equation*} 
    \gamma:[0,1]\longrightarrow Q\setminus X
    \end{equation*}
    such that    
    \begin{equation*}
    \gamma(0)\in\operatorname{relint}(F_i),\qquad \gamma(1)\in\operatorname{relint}(F_j),\qquad \gamma((0,1))\subset\operatorname{int}(Q).
    \end{equation*}
    Moreover, $\gamma$ may be chosen to meet the faces transversally.
\end{lemma}
\begin{proof}
    Let $X_Y$ denote the set of triple-junction points of $X\cap \interior(Q)$ and let $\Sigma = X_Y \cup \{v_1,v_2,v_3,v_4\}$ be the corresponding compact embedded $C^1$ one-dimensional manifold with boundary $\partial \Sigma = \{v_1,v_2,v_3,v_4\}$. 
    For each $i\in \{1,2,3,4\}$ let $V_i$ be the connected component of $Q\setminus X$ containing $\operatorname{relint} F_i$. 
    At any point $x\in X_Y$, since $X$ is modeled on the triple junction near $x$ we see that there are at most three connected components of $Q\setminus X$ which touch the point $x$. Moreover, this set of connected components is locally constant along the 1-manifold $\Sigma$. 

    Now pick a vertex $v_i$. Near $v_i$,  
     the 1-manifold $\Sigma$ consists of a simple curve $\gamma_i$. Since $\partial \Sigma = \{v_1,v_2,v_3,v_4\}$, the curve $\gamma_i$ must connect $v_i$ to another vertex $v_k$ for $i\neq k$. 
     Near $v_i$, the connected components of $Q\setminus X$ consist of the set $\{V_\ell : \ell\neq i\}$, while near $v_k$ the connected components of $Q\setminus X$ consist of the set $\{V_\ell: \ell\neq k\}$. Constancy along $\gamma_i$ implies $\{V_\ell : \ell\neq i\} = \{V_\ell: \ell\neq k\}$. It follows that $V_i$ belongs to this common set, and hence there is an index $\ell\neq i$ such that $V_i = V_\ell$. If we set $j=\ell$, then the common component $V_i =V_j$ is path-connected and so we  obtain our desired curve $\gamma$.
\end{proof}

The following proposition is the main ingredient in our proof of Theorem \ref{thm:existence}. 

\begin{proposition} \label{prop:deformation-key}
     Let $X\subset Q$ be compact with $X\cap\partial Q=T\cap\partial Q$.
Suppose there is an embedded $C^1$ curve $\gamma:[0,1]\to Q\setminus X$
such that, for some $i\neq j$,
\[
\gamma(0)\in\operatorname{relint}(F_i),\qquad
\gamma(1)\in\operatorname{relint}(F_j),\qquad
\gamma((0,1))\subset\interior(Q),
\]
and $\gamma$ meets these faces transversally. Suppose also that, for some $C>0$, 
\[
\dist(z,T\cap\partial Q)\le C\dist(z,\partial Q)
\qquad\text{for every }z\in X . 
\]    
    Then, there is a Lipschitz map $f:\R^3\to \R^3$ such that 
    \begin{equation*}
        f=\id\; \text{ on }\R^3\setminus \interior(Q),\quad f(Q) \subset Q,\quad f(X)\subset G_{ij}.
    \end{equation*}
\end{proposition}
\begin{proof}
     Since $\gamma$ is disjoint from $X$, there is a sufficiently small tubular neighborhood $N$ of $\gamma$ such that $N\cap X=\varnothing $, $N$ has the same topology of $\gamma$, and $N \cap \partial Q$ is the union of two disks centered at the endpoints of $\gamma$ and contained in the respective faces. Define
    \begin{equation*}
        M = Q \setminus \interior_Q N, \quad A = M\cap \partial Q .
    \end{equation*}
    Hence $M$ is the complement of $\interior(N)$ in $Q$ and $A$ is $\partial Q$ minus two disks. We claim that there exists a Lipschitz map $r:M \to A$ such that $r|_A={\rm id}$. We give two distinct proofs. 

    \textit{Algebraic construction.} The key observation is that $M$ is a knot complement; in particular,  $H_1(M;\Z)\cong \Z$ and $H_0(M;\Z) = \Z$; see \cite[Section 3A]{Rolfsen2003}. In our situation we can construct the isomorphisms as follows. First, $H_1(M\cap N;\Z)$ is naturally isomorphic to $H_1(A;\Z)\cong \Z$ via the two boundary circles of $M\cap N$. Meanwhile,  we have the isomorphism $H_1(M\cap N;\Z) \cong H_1(M;\Z)$ induced by the inclusion. This latter fact may be seen by the Mayer-Vietoris sequence
      \begin{equation*}
        H_2 (Q;\Z) \longrightarrow H_1(M\cap N; \Z)\longrightarrow H_1(M;\Z) \oplus H_1(N;\Z) \longrightarrow H_1(Q;\Z), 
    \end{equation*}
    and the observation that $Q$ and $N$ are topological 3-balls. Composing these two isomorphisms and the abelianization of the fundamental group, we obtain the group homomorphism \begin{equation}\label{eq:homomorphism-pi1-M-A}
        \pi_1(M)\longrightarrow  H_1(M;\Z) \longrightarrow H_1(A) \cong \pi_1(A).
    \end{equation}
    Since $A$ is a $K(\Z,1)$ and $M$ is connected, by \cite[Proposition 1B.9]{Hatcher2002} the homomorphism in \eqref{eq:homomorphism-pi1-M-A} is induced by a continuous map $r_0:M\to A$ such that $r_0|_A$ is homotopic to the identity $\id_A$.
    Applying the homotopy extension property  to the pair $(M,A)$ we obtain a continuous map $r:M\to A$ such that $r|_A = \id_A$. 
    By approximating the continuous map $r$ with Lipschitz maps equal to the identity on $A$, we may assume that $r$ is Lipschitz (see for instance \cite[Theorem 5.17]{LuukkainenVaisala1977}). 

    \textit{Geometric construction.} Up to a bi-Lipschitz map, we can identify \(A\) with \([0,1]\times \Sp^1\), using coordinates \((h,\theta)\), with the angular orientation chosen so that each circle \(\{t\}\times\Sp^1\) has linking number $+1$ with any oriented closed extension of $\gamma$. By Lemma \ref{lem: existence map} there exists \(u\in C^1(M;\Sp^1)\) with degree one on each circle \(\{t\}\times \Sp^1\subset A\). We additionally arrange that \(u|_A=\theta\). Indeed, the complex quotient
\(\theta u^{-1}|_A\) has degree zero on every circle \(\{t\}\times \Sp^1\), and therefore admits a
Lipschitz lift \(\psi:A\to\mathbb R\), such that $\theta=e^{ i\psi}u$ on $A$. By McShane's extension theorem, extend \(\psi\) to a
Lipschitz function \(\Psi:M\to\mathbb R\).
Replacing \(u\) by \(e^{ i\Psi}u\), we obtain
a Lipschitz map satisfying \(u\equiv \theta\) on $A$.

Again, by McShane's theorem followed by a truncation to
\([0,1]\), extend the height function \(h|_A\) to a Lipschitz function
\(\widetilde h:M\to[0,1]\). Under the above identification, the map $r:M\to A$ defined by $r(x)=(\widetilde h(x),u(x))$ is Lipschitz and satisfies
\(r(a)=(h(a),\theta(a))=a\) for every \(a\in A\). Thus \(r\) is the desired Lipschitz retraction.

\begin{figure}[H]
    \centering
    \resizebox{0.56\textwidth}{!}{\input{retraction_map.tikz}}
    \caption{Construction of the retraction $r:M\to A$ from the map $u$.}
    \label{fig: retraction with u}
\end{figure}
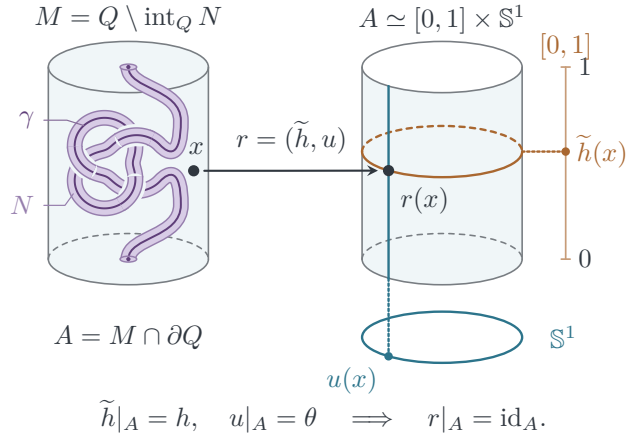

    Once we have the existence of such $r:M \to A$, we can proceed as in \cite[p. 211]{Liang2014} to construct the desired map $f$. 
    To be precise, let $p:A\to G_{ij}$ be the Lipschitz map which projects radially away from the centers of the two disks of $\partial Q \cap N$.
    Since $X\subset M$, we can define the Lipschitz map $g=(p\circ r)|_X$ which fixes $T\cap \partial Q$ and maps $X$ into $G_{ij}$. 
    We can then extend $g$ to a Lipschitz map $\overline g: X\cup \partial Q\to \partial Q$ by setting 
    $\overline g|_X= g$ and $\overline g|_{\partial Q} =\id_{\partial Q}$. 
    Since $g|_{T\cap \partial Q}=\id_{T\cap \partial Q}$ the map is well-defined, and one can check that the condition $\dist(z,T\cap \partial Q) \leq C  \dist(z,\partial Q)$ for $z\in X$ ensures that $\overline g$ is Lipschitz. 

    We can now apply McShane's theorem to $\overline g$ to obtain a Lipschitz function $\R^3\to \R^3$ agreeing with $\overline g$ on $X\cup \partial Q$. Composing with the nearest-point projection to the convex set $Q$ and extending by the identity outside of $Q$, we obtain the desired Lipschitz map $f$.
\end{proof}

The main topological input to this construction is the following fact.
\begin{lemma}\label{lem: existence map}
   Let $M,Q$, and $\gamma$ be as in the proof of Proposition \ref{prop:deformation-key}. Then there is a map $u\in C^1(M;\Sp^1)$ such that, for every closed loop $\ell : \Sp^1 \to M$, we have
    \begin{equation*}
        \deg (u\circ \ell ) = {\rm link}(\ell , \widehat \gamma) , 
    \end{equation*}
    where $\widehat \gamma $ is any oriented closed extension of $\gamma$ outside $Q$ with
$\widehat\gamma\cap Q=\gamma([0,1])$. 
\end{lemma}
\begin{proof}
    For the closed curve $\widehat\gamma$, it is classical that there is a smooth map
$ u:\R^3\setminus\widehat\gamma\to\Sp^1$ whose degree along
each loop equals its linking number with $\widehat\gamma$; see for example \cite[Theorem 3.3]{BM21} or \cite[Theorem 4.4]{ABO1}. The restriction $u|_M$ has the desired property. 
\end{proof}

We may now prove Theorem~\ref{thm:existence}.

\begin{proof}[Proof of Theorem~\ref{thm:existence}]
    If the theorem is false, after translating, dilating, and rotating we obtain   a sequence of reduced two-dimensional Almgren-minimizing sets $E_k$ in $U_k\supset \overline{B_8 }$ such that $d_{0,8}(E_k,T) \to 0$ and with $E_k\cap B_2$ containing no tetrahedral point. 

    Apply Lemma~\ref{lem:annular-diffeo} to $E_k\cap B_4$ and set
$X_k:=\Phi_k^{-1}(E_k)\cap Q$ and $Q_k:=\Phi_k(Q)$.
For all sufficiently large $k$, we have that $Q_k\subset B_2$. Since $E_k\cap B_2$ contains no $T$ point, Taylor's theorem implies
that every point of $E_k\cap Q_k$ is locally diffeomorphic to a plane
or a triple junction. 

Thus $X_k$ satisfies the hypotheses of
Lemma~\ref{lem:topological-path}. Applying that lemma and then
Proposition~\ref{prop:deformation-key}, we obtain Lipschitz maps
$f_k:\R^3\to\R^3$ and distinct indices
$i_k,j_k\in\{1,2,3,4\}$ satisfying
    \begin{equation*}
        f_k = \id \text{ on } \R^3 \setminus \interior(Q), \quad f_k(Q)\subset Q,\quad f_k(X_k)\subset G_{i_k,j_k}.
    \end{equation*}
    If we set $\varphi_k = \Phi_k \circ f_k \circ \Phi_k^{-1}$, then $\varphi_k(Q_k) = Q_k$ and  $\varphi_k(E_k)$ is an Almgren competitor for $E_k$ agreeing with $E_k$ outside of $Q_k$.
    Moreover by the construction of $f_k$, 
    \begin{equation*}
        \varphi_k(E_k) \cap Q_k = \Phi_k(f_k(X_k))\subset\Phi_k(G_{i_k,j_k}) . 
    \end{equation*}
Applying minimality, we obtain
    \begin{align*}
        \cH^2 (E_k\cap Q_k)\leq \cH^2(\varphi_k(E_k)\cap Q_k) & \leq \cH^2(\Phi_k(G_{i_k,j_k})) \\ & \leq \left( 1+ \|D\Phi_k-I\|_{L^\infty}\right)^2 \cH^2(G_{i_kj_k}).
    \end{align*}
  By Lemma~\ref{lem:annular-diffeo},
$\cH^2(E_k\cap Q_k)\to2\sqrt2$ and
$\|D\Phi_k-I\|_{L^\infty}\to0$. Since
$\cH^2(G_{i_kj_k})=4\sqrt3/3$ for every $k$, passing to the limit gives
\[
2\sqrt2\le\frac{4\sqrt3}{3}<2\sqrt2,
\]
contradiction.  
\end{proof}

\appendix
\setcounter{section}{1}
\setcounter{equation}{0}
\section*{Appendix}

\setcounter{theorem}{0}
\renewcommand{\thetheorem}{\Alph{theorem}}
\renewcommand{\theHtheorem}{appendix.\arabic{theorem}}

The following is the partial result on the uniqueness of $T$ with prescribed boundary datum that can be deduced from our Proposition~\ref{prop:deformation-key}.

\begin{theorem}\label{thm: prescribed-boundary}
Let $E\subset Q$ be a reduced, two-dimensional Almgren-minimal set in $\interior(Q)$ with $E\cap\partial Q= T \cap \partial Q =: \Gamma$ and $\mathcal{H}^2(E)<+\infty$. Assume that there exists a $C^1$
diffeomorphism $\Psi:\R^3\to\R^3$ and an open neighborhood $U$
of $\partial Q$ such that
\begin{equation}\label{eq: boundary hp}
    \Psi(Q)=Q,\qquad \Psi(\Gamma)=\Gamma,\qquad
    E\cap U=\Psi(T\cap Q)\cap U.
\end{equation}
Then $E \equiv T\cap Q$.
\end{theorem}

\begin{proof}
Taylor's theorem and the boundary hypothesis give a
piecewise $C^1$ surface, smooth on each open sheet.
Let $\nu$ denote its outward unit conormal along the six open edges
of $\Gamma$. Interior stationarity, the divergence theorem on $E$, and cancellation of the
conormals along the singular curves give
\begin{equation*}
    \int_E \operatorname{div}_E Z\,d\cH^2
    =\int_{\Gamma} Z\cdot\nu\,d\cH^1
    \qquad\text{for every } Z\in C_c^1(\R^3;\R^3).
\end{equation*}
In particular, $E$ is strongly stationary with respect to $\Gamma$
in the sense of \cite{GulliverYamada2006}. Taking $Z(x)=x$ near $Q$
yields
\[
    2\cH^2(E)=\int_\Gamma x\cdot\nu\,d\cH^1.
\]

Let $m$ be the midpoint of an edge of $Q$. The normal section of
$Q$ along this edge is a wedge with inward bisector $-m$ and
half-angle
\[
    \beta=\frac12\arccos(1/3),\qquad
    |m|=1/ \sqrt3 ,\qquad \cos\beta=\sqrt{2/3}.
\]
The boundary hypothesis \eqref{eq: boundary hp} implies that $-\nu$ lies strictly inside this wedge. Since $x-m$ is tangent to the edge for $x\in \Gamma$, we get
\[
    x\cdot\nu=m\cdot\nu>|m|\cos\beta=\sqrt2 /3 .
\]
Each edge has length $\sqrt{8/3}$, and hence
\begin{equation}\label{eq:boundary-area-gap}
    \cH^2(E)>
    \frac12 \cdot 6\sqrt{\frac83} \cdot \frac{\sqrt2}{3}
    =\frac{4\sqrt3}{3}
    =\cH^2(G_{ij}).
\end{equation}

Suppose now that $E$ has no interior $T$ point. The set
$X=\Psi^{-1}(E)$ satisfies the hypotheses of
Lemma~\ref{lem:topological-path}. Since $\Psi$ permutes the faces
of $Q$, the image of the curve supplied by that lemma joins two
distinct faces in $Q\setminus E$. Moreover, the boundary hypothesis \eqref{eq: boundary hp} gives
\[
    \dist(z,\Gamma)\le C\dist(z,\partial Q)
    \qquad\text{for every }z\in E.
\]
Proposition~\ref{prop:deformation-key} therefore gives a
Lipschitz map $f:\R^3\to\R^3$ satisfying
\[
    f=\id\text{ on }\R^3\setminus\interior(Q),\qquad
    f(Q)\subset Q,\qquad f(E)\subset G_{ij}.
\]
We approximate $f$ by admissible deformations in the interior of $Q$. Set
$d(x)=\dist(x,\partial Q)$, and choose a Lipschitz function
$\chi:[0,\infty)\to[0,1]$ equal to $0$ on $[0,1]$ and to $1$
on $[2,\infty)$. For $0<\delta<1$, set $ \eta_\delta(x):=\chi(d(x)/\delta)$, define
\[
    f_\delta(x)=(1-\eta_\delta(x))x
        +\eta_\delta(x)(1-\delta)f(x),
\]
and extend $f_\delta$ by the identity outside $Q$.
This map fixes a neighborhood of $\partial Q$ and maps $\interior(Q)$ into
itself. Its straight-line homotopy is therefore admissible.
Since $f=\id$ on $\partial Q$,
\[
    |f(x)-x|\le(\Lip(f)+1)d(x).
\]
Thus $|(1-\delta)f(x)-x|\le C\delta$ on the transition layer
$\{\delta\le d\le2\delta\}$, and consequently
$\Lip(f_\delta)\le C_1 $ independently of $\delta$.
Minimality and the area estimate for Lipschitz maps give
\[
    \cH^2(E)
    \le\cH^2(f_\delta(E))
    \le(1-\delta)^2\cH^2(G_{ij})
       +C_1^2\cH^2(E\cap\{d\le2\delta\}).
\]
The last term tends to zero because $\cH^2(E)<\infty$ and
$\cH^2(E\cap\partial Q)=0$. Letting $\delta\downarrow0$
contradicts the mass lower-bound \eqref{eq:boundary-area-gap}. Hence there exists
a $T$ point $p\in E\cap\interior(Q)$.

Let $C_p(\Gamma)$ be the cone from $p$ over $\Gamma$. By \cite[Corollary 1 and Example 2]{GulliverYamada2006} we get
\[
    \theta_T=\theta_E(p)\le \frac{1}{2} \mathcal{C}_{\rm tot}(\Gamma)= \theta_T , 
\]
where $\mathcal{C}_{\rm tot}(\Gamma)$ is the total curvature of $\Gamma$ as in \cite[eq. (2)]{GulliverYamada2006}. By the equality case of
\cite[Theorem 1]{GulliverYamada2006}, $E=C_p(\Gamma)$ is the cone over $\Gamma$ from $p$.

Since $p$ is a $T$ point, the unit vectors
\[
    u_i=\frac{v_i-p}{|v_i-p|} \in \Sp^2 ,\qquad r_i=|v_i-p|>0,
\]
satisfy $u_i\cdot u_j=-1/3$ for $i\ne j$. All edges of $Q$
have the same length, so, for distinct $i,j,k$,
\[
    0=|v_i-v_j|^2-|v_i-v_k|^2
     =(r_j-r_k)\left(r_j+r_k+\frac23r_i\right), 
\]
and hence $r_j=r_k$. This clearly implies $p=0$, thus $E=C_0(\Gamma)=T\cap Q$.
\end{proof}

\subsubsection*{Acknowledgments}

M.C. acknowledges support from the Anthony and Sharon Lee Foundation and is very grateful to Princeton University and the Institute for Advanced Study for their kind hospitality during the realization of this work.

Z.L. acknowledges the support of the National Science Foundation Graduate Research Fellowship Program under Grant No. DGE-2444107. Any opinions, findings, and conclusions or recommendations expressed in this material are those of the author(s) and do not necessarily reflect the views of the National Science Foundation.

\bibliographystyle{alphaurl-initials}

\bibliography{references}
\end{document}

%% file: y_mesh_support.tex
\makeatletter
\@ifundefined{YMeshGraphic}{%
  \newbox\YMeshBox
  \newcommand{\YMeshGraphic}[1]{%
    \begingroup
    \pdfcompresslevel=0
    \immediate\pdfobj stream attr{%
      /ShadingType 4 /ColorSpace /DeviceRGB
      /BitsPerCoordinate 24 /BitsPerComponent 8 /BitsPerFlag 8
      /Decode [0 600 0 600 0 1 0 1 0 1]
      /Filter [/ASCII85Decode /FlateDecode]
    } file{#1}%
    \edef\YShadingObject{\the\pdflastobj}%
    \setbox\YMeshBox=\hbox to 600bp{%
      \vrule width0pt height600bp depth0pt
      \pdfliteral direct{/YShading sh}\hss
    }%
    \pdfcompresslevel=9
    \immediate\pdfxform resources{%
      /Shading << /YShading \YShadingObject\space 0 R >>
    }\YMeshBox
    \pdfrefxform\pdflastxform
    \endgroup
  }%
}{}
\makeatother
\definecolor{filmblue}{RGB}{120,188,242}
\definecolor{filmorange}{RGB}{204,128,72}
\definecolor{filmgreen}{RGB}{124,221,147}
\definecolor{filmlabel}{RGB}{55,62,68}

%% file: y_set.tikz
\begingroup%
\selectcolormodel{rgb}%
\definecolor{YFilmWire}{RGB}{185,185,183}%
\definecolor{YFilmWhite}{RGB}{245,245,242}%
\definecolor{YFilmJunction}{RGB}{201,200,194}%
%
\def\YFilmAdd#1#2#3{%
  dup #3 mul 3 -1 roll add exch
  dup #2 mul 4 -1 roll add 3 1 roll
  #1 mul 4 -1 roll add 3 1 roll
}%
\def\YFilmGlow#1#2#3#4#5#6#7{%
  4 index #1 sub #3 mul dup mul
  4 index #2 sub #4 mul dup mul
  add neg 2.718281828 exch exp
  \YFilmAdd{#5}{#6}{#7}
}%
\pgfdeclarefunctionalshading{YFilmBack}{\pgfpoint{-50bp}{-50bp}}{\pgfpoint{50bp}{50bp}}{}{
  25 exch sub 50 div exch 25 add 50 div exch
  0.400000 0.603000 0.781000
  4 index \YFilmAdd{0.045000}{0.062000}{0.074000}
  3 index \YFilmAdd{0.015000}{0.020000}{0.030000}
  \YFilmGlow{0.640000}{0.260000}{2.300000}{2.300000}{0.072000}{0.110000}{0.122000}
  \YFilmGlow{0.060000}{0.700000}{4.000000}{1.400000}{-0.045000}{-0.068000}{-0.090000}
  \YFilmGlow{0.680000}{0.320000}{9.000000}{5.000000}{0.045000}{0.048000}{0.035000}
  5 -2 roll pop pop
}%
\pgfdeclarefunctionalshading{YFilmLeft}{\pgfpoint{-50bp}{-50bp}}{\pgfpoint{50bp}{50bp}}{}{
  25 exch sub 50 div exch 25 add 50 div exch
  0.480000 0.714000 0.929000
  4 index \YFilmAdd{-0.080000}{-0.120000}{-0.140000}
  3 index \YFilmAdd{0.035000}{0.040000}{0.034000}
  \YFilmGlow{0.300000}{0.720000}{2.400000}{2.200000}{0.060000}{0.090000}{0.064000}
  \YFilmGlow{0.400000}{0.400000}{5.000000}{4.000000}{0.024000}{0.040000}{0.030000}
  \YFilmGlow{0.970000}{0.400000}{6.000000}{1.400000}{-0.055000}{-0.080000}{-0.105000}
  \YFilmGlow{0.280000}{0.660000}{9.000000}{5.000000}{0.025000}{0.032000}{0.009000}
  5 -2 roll pop pop
}%
\pgfdeclarefunctionalshading{YFilmFront}{\pgfpoint{-50bp}{-50bp}}{\pgfpoint{50bp}{50bp}}{}{
  25 exch sub 50 div exch 25 add 50 div exch
  0.384000 0.578000 0.745000
  4 index \YFilmAdd{0.055000}{0.067000}{0.077000}
  3 index \YFilmAdd{0.020000}{0.025000}{0.020000}
  \YFilmGlow{0.750000}{0.620000}{2.000000}{1.600000}{0.035000}{0.055000}{0.075000}
  \YFilmGlow{0.680000}{0.570000}{8.000000}{4.000000}{0.045000}{0.060000}{0.063000}
  \YFilmGlow{0.050000}{0.380000}{5.000000}{1.800000}{-0.025000}{-0.034000}{-0.038000}
  5 -2 roll pop pop
}%
%
\begin{tikzpicture}[
    x=0.020161290323cm,y=-0.020161290323cm,
    baseline=(current bounding box.center),
    line cap=round,line join=round,
    frame/.style={draw=YFilmWire,line width=1.12pt},
    frame light/.style={draw=YFilmWhite,line width=0.68pt},
    junction/.style={draw=YFilmJunction,line width=1.05pt},
    junction light/.style={draw=YFilmWhite,line width=0.55pt}
  ]
  \path[use as bounding box] (-3cm,-2.7cm) rectangle (3cm,2.7cm);
  \fill[white] (-3cm,-2.7cm) rectangle (3cm,2.7cm);
  \begin{scope}[shift={(-151.00000000,-160.00000000)}]

  \coordinate (Jtop) at (157,87);
  \coordinate (Jbottom) at (157,232);
  \coordinate (Ltop) at (75,59);
  \coordinate (Lbottom) at (75,204);
  \coordinate (Rtop) at (227,36);
  \coordinate (Rbottom) at (227,181);
  \coordinate (Ftop) at (191,139);
  \coordinate (Fbottom) at (191,284);

  \begin{scope}
    \clip (Jtop)--(Rtop)--(Rbottom)--(Jbottom) -- cycle;
    \shade[shading=YFilmBack] (157.0000,36.0000) rectangle (227.0000,232.0000);
  \end{scope}
  \draw[frame] (Jtop)--(Rtop)--(Rbottom)--(Jbottom);
  \draw[frame light] (Jtop)--(Rtop)--(Rbottom)--(Jbottom);

  \begin{scope}
    \clip (Jtop)--(Ltop)--(Lbottom)--(Jbottom) -- cycle;
    \shade[shading=YFilmLeft] (75.0000,59.0000) rectangle (157.0000,232.0000);
  \end{scope}
  \draw[frame] (Jtop)--(Ltop)--(Lbottom)--(Jbottom);
  \draw[frame light] (Jtop)--(Ltop)--(Lbottom)--(Jbottom);

  \begin{scope}
    \clip (Jtop)--(Ftop)--(Fbottom)--(Jbottom) -- cycle;
    \shade[shading=YFilmFront] (157.0000,87.0000) rectangle (191.0000,284.0000);
  \end{scope}
  \draw[frame] (Jtop)--(Ftop)--(Fbottom)--(Jbottom);
  \draw[frame light] (Jtop)--(Ftop)--(Fbottom)--(Jbottom);

  \draw[junction] (Jtop)--(Jbottom);
  \draw[junction light] (Jtop)--(Jbottom);
  \end{scope}
\end{tikzpicture}%
\endgroup%

%% file: t_set.tikz
\begingroup%
\selectcolormodel{rgb}%
\definecolor{TFilmWire}{RGB}{185,185,183}%
\definecolor{TFilmWhite}{RGB}{245,245,242}%
\definecolor{TFilmJunction}{RGB}{201,200,194}%
%
\def\TFilmAdd#1#2#3{%
  dup #3 mul 3 -1 roll add exch
  dup #2 mul 4 -1 roll add 3 1 roll
  #1 mul 4 -1 roll add 3 1 roll
}%
\def\TFilmGlow#1#2#3#4#5#6#7{%
  4 index #1 sub #3 mul dup mul
  4 index #2 sub #4 mul dup mul
  add neg 2.718281828 exch exp
  \TFilmAdd{#5}{#6}{#7}
}%
\pgfdeclarefunctionalshading{TFilmDeep}{\pgfpoint{-50bp}{-50bp}}{\pgfpoint{50bp}{50bp}}{}{
  25 exch sub 50 div exch 25 add 50 div exch
  0.410000 0.620000 0.805000
  4 index \TFilmAdd{0.018000}{0.022000}{0.030000}
  3 index \TFilmAdd{0.018000}{0.028000}{0.036000}
  \TFilmGlow{0.500000}{0.600000}{2.000000}{2.000000}{0.040000}{0.055000}{0.070000}
  5 -2 roll pop pop
}%
\pgfdeclarefunctionalshading{TFilmBackLeft}{\pgfpoint{-50bp}{-50bp}}{\pgfpoint{50bp}{50bp}}{}{
  25 exch sub 50 div exch 25 add 50 div exch
  0.400000 0.606000 0.791000
  4 index \TFilmAdd{0.018000}{0.033000}{0.050000}
  3 index \TFilmAdd{0.046000}{0.067000}{0.075000}
  \TFilmGlow{0.280000}{0.570000}{2.300000}{2.000000}{0.035000}{0.055000}{0.065000}
  \TFilmGlow{0.870000}{0.630000}{3.200000}{2.600000}{-0.068000}{-0.097000}{-0.119000}
  \TFilmGlow{0.320000}{0.570000}{9.000000}{5.000000}{0.024000}{0.036000}{0.032000}
  5 -2 roll pop pop
}%
\pgfdeclarefunctionalshading{TFilmRight}{\pgfpoint{-50bp}{-50bp}}{\pgfpoint{50bp}{50bp}}{}{
  25 exch sub 50 div exch 25 add 50 div exch
  0.446000 0.673000 0.882000
  4 index \TFilmAdd{0.022000}{0.028000}{0.025000}
  3 index \TFilmAdd{0.040000}{0.050000}{0.045000}
  \TFilmGlow{0.570000}{0.380000}{2.100000}{2.100000}{0.055000}{0.083000}{0.060000}
  \TFilmGlow{0.080000}{0.770000}{4.000000}{3.000000}{-0.100000}{-0.140000}{-0.180000}
  \TFilmGlow{0.600000}{0.450000}{9.000000}{5.000000}{0.040000}{0.047000}{0.023000}
  5 -2 roll pop pop
}%
\pgfdeclarefunctionalshading{TFilmBottomLeft}{\pgfpoint{-50bp}{-50bp}}{\pgfpoint{50bp}{50bp}}{}{
  25 exch sub 50 div exch 25 add 50 div exch
  0.520000 0.750000 0.948000
  4 index \TFilmAdd{-0.040000}{-0.050000}{-0.030000}
  3 index \TFilmAdd{0.028000}{0.030000}{0.012000}
  \TFilmGlow{0.350000}{0.450000}{2.000000}{2.500000}{0.020000}{0.040000}{0.025000}
  \TFilmGlow{0.960000}{0.000000}{4.000000}{3.000000}{-0.038000}{-0.056000}{-0.064000}
  5 -2 roll pop pop
}%
\pgfdeclarefunctionalshading{TFilmBottomRight}{\pgfpoint{-50bp}{-50bp}}{\pgfpoint{50bp}{50bp}}{}{
  25 exch sub 50 div exch 25 add 50 div exch
  0.465000 0.710000 0.936000
  4 index \TFilmAdd{0.035000}{0.015000}{0.003000}
  3 index \TFilmAdd{0.015000}{0.020000}{0.006000}
  \TFilmGlow{0.500000}{0.550000}{2.100000}{2.000000}{0.050000}{0.065000}{0.025000}
  \TFilmGlow{0.280000}{0.050000}{3.000000}{4.000000}{-0.030000}{-0.045000}{-0.060000}
  \TFilmGlow{0.550000}{0.550000}{9.000000}{4.000000}{0.033000}{0.041000}{0.010000}
  5 -2 roll pop pop
}%
\pgfdeclarefunctionalshading{TFilmFront}{\pgfpoint{-50bp}{-50bp}}{\pgfpoint{50bp}{50bp}}{}{
  25 exch sub 50 div exch 25 add 50 div exch
  0.361000 0.542000 0.701000
  4 index \TFilmAdd{0.050000}{0.060000}{0.074000}
  3 index \TFilmAdd{0.020000}{0.025000}{0.033000}
  \TFilmGlow{0.700000}{0.500000}{2.400000}{1.600000}{0.065000}{0.092000}{0.121000}
  \TFilmGlow{0.100000}{0.300000}{4.000000}{2.500000}{-0.021000}{-0.030000}{-0.038000}
  \TFilmGlow{0.690000}{0.520000}{8.000000}{5.000000}{0.024000}{0.037000}{0.037000}
  5 -2 roll pop pop
}%
%
\begin{tikzpicture}[
    x=0.024969181787cm,y=-0.024969181787cm,
    baseline=(current bounding box.center),
    line cap=round,line join=round,
    frame/.style={draw=TFilmWire,line width=1.12pt},
    frame light/.style={draw=TFilmWhite,line width=0.68pt},
    junction/.style={draw=TFilmJunction,line width=1.05pt},
    junction light/.style={draw=TFilmWhite,line width=0.55pt}
  ]
  \path[use as bounding box] (-3cm,-2.7cm) rectangle (3cm,2.7cm);
  \fill[white] (-3cm,-2.7cm) rectangle (3cm,2.7cm);
  \begin{scope}[shift={(-633.96916000,-193.96799500)}]
  \coordinate (A) at (611.03203,93.84457);
  \coordinate (B) at (530.09415,250.04894);
  \coordinate (C) at (576.02966,294.09142);
  \coordinate (D) at (737.84417,260.01507);
  \coordinate (O) at (613.75000,224.50000);

  \begin{scope}
    \clip (O)--(B)--(D) -- cycle;
    \shade[shading=TFilmDeep] (530.0941,224.5000) rectangle (737.8442,260.0151);
  \end{scope}

  \begin{scope}
    \clip (O)--(A)--(B) -- cycle;
    \shade[shading=TFilmBackLeft] (530.0941,93.8446) rectangle (613.7500,250.0489);
  \end{scope}
  \begin{scope}
    \clip (O)--(A)--(D) -- cycle;
    \shade[shading=TFilmRight] (611.0320,93.8446) rectangle (737.8442,260.0151);
  \end{scope}

  \begin{scope}
    \clip (O)--(B)--(C) -- cycle;
    \shade[shading=TFilmBottomLeft] (530.0941,224.5000) rectangle (613.7500,294.0914);
  \end{scope}

  \draw[junction] (B)--(O);
  \draw[junction light] (B)--(O);

  \begin{scope}
    \clip (O)--(C)--(D) -- cycle;
    \shade[shading=TFilmBottomRight] (576.0297,224.5000) rectangle (737.8442,294.0914);
  \end{scope}

  \begin{scope}
    \clip (O)--(A)--(C) -- cycle;
    \shade[shading=TFilmFront] (576.0297,93.8446) rectangle (613.7500,294.0914);
  \end{scope}

  \draw[junction] (A)--(O)--(C);
  \draw[junction light] (A)--(O)--(C);
  \draw[junction] (O)--(D);
  \draw[junction light] (O)--(D);

  \draw[frame] (A)--(B)--(C)--(D)--cycle;
  \draw[frame light] (A)--(B)--(C)--(D)--cycle;
  \draw[frame] (A)--(C);
  \draw[frame light] (A)--(C);
  \end{scope}
\end{tikzpicture}%
\endgroup%

%% file: retraction_map.tikz
\begingroup
\definecolor{retrink}{RGB}{48,55,63}
\definecolor{retrheight}{RGB}{169,103,48}
\definecolor{retrcircle}{RGB}{44,116,139}
\definecolor{retrsurface}{RGB}{240,246,248}
\begin{tikzpicture}[
  x=0.85cm,y=0.85cm,
  font=\footnotesize,
  text=retrink,
  line cap=round,line join=round,
  >=stealth,
  boundary/.style={draw=retrink!65,line width=0.6pt},
  hidden/.style={dash pattern=on 2pt off 2pt},
  height/.style={draw=retrheight,line width=0.9pt},
  circle coordinate/.style={draw=retrcircle,line width=0.9pt},
  point/.style={circle,fill=retrink,inner sep=1.5pt,outer sep=0pt}
]
  \def\cx{5.10}
  \def\rad{1.25}
  \def\dep{0.40}
  \def\base{0.80}
  \def\top{3.80}
  \def\level{2.48}
  \def\phase{228}
  \pgfmathsetmacro{\ux}{\cx+\rad*cos(\phase)}
  \pgfmathsetmacro{\uy}{\dep*sin(\phase)}
  \coordinate (source) at (1.22,2.18);
  \coordinate (image) at (\ux,\level+\uy);
  \coordinate (angle) at (\ux,-0.42+\uy);

\definecolor{retrtube}{RGB}{216,204,230}
\definecolor{retrtubeedge}{RGB}{150,126,173}
\definecolor{retrknot}{RGB}{95,62,121}

\path[fill=retrsurface]
  (-1.05,0.80)
  arc[start angle=180,end angle=360,x radius=1.25,y radius=0.40]
  -- (1.45,3.80)
  arc[start angle=0,end angle=-180,x radius=1.25,y radius=0.40]
  -- cycle;
\path[fill=retrsurface] (0.20,3.80)
  ellipse[x radius=1.25,y radius=0.40];
\draw[boundary,hidden] (-1.05,0.80)
  arc[start angle=180,end angle=0,x radius=1.25,y radius=0.40];
\draw[boundary] (-1.05,0.80) -- (-1.05,3.80);
\draw[boundary] (1.45,0.80) -- (1.45,3.80);
\draw[boundary] (0.20,3.80)
  ellipse[x radius=1.25,y radius=0.40];
\node at (0.20,4.55) {$M=Q\setminus\operatorname{int}_Q N$};

\pgfmathsetmacro{\retrstartx}{0.27*(2+cos(36))*cos(24)}
\pgfmathsetmacro{\retrstarty}{2.30+0.25*(2+cos(36))*sin(24)}
\pgfmathsetmacro{\retrendx}{\retrstartx}
\pgfmathsetmacro{\retrendy}{2.30-0.25*(2+cos(36))*sin(24)}
\def\retrArcPath{
  (0.20,3.80)
  .. controls (0.20,3.45) and (1.02,3.47) .. (1.02,2.96)
  .. controls (1.02,2.38) and (0.84,2.43)
     .. (\retrstartx,\retrstarty)
  plot[domain=12:348,samples=337,variable=\t]
    ({0.27*(2+cos(3*\t))*cos(2*\t)},
     {2.30+0.25*(2+cos(3*\t))*sin(2*\t)})
  .. controls (0.84,2.17) and (1.02,2.22) .. (1.02,1.65)
  .. controls (1.02,1.13) and (0.20,1.15) .. (0.20,0.80)
}
\draw[draw=retrtubeedge,line width=6.0pt,line cap=butt] \retrArcPath;
\draw[draw=retrtube,line width=4.9pt,line cap=butt] \retrArcPath;
\draw[draw=retrknot,line width=0.75pt,line cap=butt] \retrArcPath;

\foreach \retrcross in {30,150,270} {
  \pgfmathsetmacro{\retrfirst}{\retrcross-12}
  \pgfmathsetmacro{\retrlast}{\retrcross+12}
  \pgfmathsetmacro{\retrmaskfirst}{\retrcross-5}
  \pgfmathsetmacro{\retrmasklast}{\retrcross+5}
  \draw[draw=retrsurface,line width=8.4pt]
    plot[domain=\retrmaskfirst:\retrmasklast,samples=11,variable=\t]
      ({0.27*(2+cos(3*\t))*cos(2*\t)},
       {2.30+0.25*(2+cos(3*\t))*sin(2*\t)});
  \draw[draw=retrtubeedge,line width=6.0pt,line cap=butt]
    plot[domain=\retrfirst:\retrlast,samples=25,variable=\t]
      ({0.27*(2+cos(3*\t))*cos(2*\t)},
       {2.30+0.25*(2+cos(3*\t))*sin(2*\t)});
  \draw[draw=retrtube,line width=4.9pt,line cap=butt]
    plot[domain=\retrfirst:\retrlast,samples=25,variable=\t]
      ({0.27*(2+cos(3*\t))*cos(2*\t)},
       {2.30+0.25*(2+cos(3*\t))*sin(2*\t)});
  \draw[draw=retrknot,line width=0.75pt,line cap=butt]
    plot[domain=\retrfirst:\retrlast,samples=25,variable=\t]
      ({0.27*(2+cos(3*\t))*cos(2*\t)},
       {2.30+0.25*(2+cos(3*\t))*sin(2*\t)});
}

\foreach \retrendy in {0.80,3.80} {
  \path[draw=retrtubeedge,fill=retrtube,line width=0.55pt]
    (0.20,\retrendy) ellipse[x radius=0.13,y radius=0.046];
  \fill[retrknot] (0.20,\retrendy) circle[radius=0.8pt];
}
\draw[boundary] (-1.05,0.80)
  arc[start angle=180,end angle=360,x radius=1.25,y radius=0.40];

\node[anchor=east,text=retrknot,inner sep=2pt]
  at (-1.18,2.97) {$\gamma$};
\draw[draw=retrknot!75,line width=0.45pt]
  (-1.12,2.97) -- (-0.76,2.97) -- (-0.496,2.887);
\node[anchor=east,text=retrtubeedge,inner sep=2pt]
  at (-1.18,1.64) {$N$};
\draw[draw=retrtubeedge,line width=0.45pt]
  (-1.12,1.64) -- (-0.79,1.64) -- (-0.66,1.82);
\node at (0.20,-0.42) {$A=M\cap\partial Q$};

  \path[fill=retrsurface]
    (\cx-\rad,\base)
    arc[start angle=180,end angle=360,x radius=\rad,y radius=\dep]
    -- (\cx+\rad,\top)
    arc[start angle=0,end angle=-180,x radius=\rad,y radius=\dep]
    -- cycle;
  \draw[boundary,hidden]
    (\cx-\rad,\base)
    arc[start angle=180,end angle=0,x radius=\rad,y radius=\dep];
  \draw[boundary] (\cx-\rad,\base) -- (\cx-\rad,\top);
  \draw[boundary] (\cx+\rad,\base) -- (\cx+\rad,\top);
  \draw[boundary] (\cx,\top) ellipse[x radius=\rad,y radius=\dep];
  \draw[boundary]
    (\cx-\rad,\base)
    arc[start angle=180,end angle=360,x radius=\rad,y radius=\dep];
  \node at (\cx,4.55) {$A\simeq[0,1]\times\mathbb S^1$};

  \draw[height,hidden]
    (\cx-\rad,\level)
    arc[start angle=180,end angle=0,x radius=\rad,y radius=\dep];
  \draw[height]
    (\cx-\rad,\level)
    arc[start angle=180,end angle=360,x radius=\rad,y radius=\dep];

  \draw[draw=retrheight!65,line width=0.65pt] (7.03,\base) -- (7.03,\top);
  \foreach \yy/\lab in {0.80/0,3.80/1} {
    \draw[draw=retrheight!65,line width=0.65pt] (6.96,\yy) -- (7.10,\yy);
    \node[anchor=west,inner sep=3pt] at (7.10,\yy) {$\lab$};
  }
  \draw[height,densely dotted] (\cx+\rad,\level) -- (7.03,\level);
  \fill[retrheight] (7.03,\level) circle[radius=1.6pt];
  \node[anchor=west,text=retrheight,inner sep=4pt] at (7.03,\level)
    {$\widetilde h(x)$};
  \node[text=retrheight] at (7.03,4.12) {$[0,1]$};

  \draw[circle coordinate] (\ux,\base+\uy) -- (\ux,\top+\uy);
  \draw[circle coordinate,densely dotted] (angle) -- (\ux,\base+\uy);
  \draw[circle coordinate]
    (\cx,-0.42) ellipse[x radius=\rad,y radius=\dep];
  \fill[retrcircle] (angle) circle[radius=1.6pt];
  \node[anchor=north east,text=retrcircle,inner sep=4pt] at (angle) {$u(x)$};
  \node[anchor=west,text=retrcircle] at (6.62,-0.42) {$\mathbb S^1$};

  \node[point] at (source) {};
  \node[anchor=south,inner sep=6pt,xshift=1pt] at (source) {$x$};
  \draw[->,draw=retrink,line width=0.8pt,shorten <=4pt,shorten >=4pt]
    (source) .. controls (2.10,2.18) and (3.10,2.18) .. (image);
  \node[inner sep=3pt] at (2.75,2.71) {$r=(\widetilde h,u)$};
  \node[point] at (image) {};
  \node[anchor=north west,inner sep=0pt]
    at (\ux+0.18,\level+\uy-0.25) {$r(x)$};

  \node at (3.25,-1.65)
    {$\widetilde h|_A=h,\quad u|_A=\theta
      \quad\Longrightarrow\quad r|_A=\mathrm{id}_A.$};
\end{tikzpicture}
\endgroup